\documentclass[12pt,reqno,fleqn,a4paper]{article}
\usepackage{amsmath,amssymb,amsfonts,amsthm,enumitem,cite}
\usepackage{setspace,xcolor}
\usepackage{graphicx}
\usepackage{multicol}
\usepackage{lipsum,lineno}
\usepackage{ragged2e}
\usepackage{multirow}
\usepackage{bm}
\usepackage[margin=1.0in]{geometry}
\usepackage[utf8]{inputenc}
\usepackage[T1]{fontenc}
\usepackage{mathtools}
\usepackage[font=small,labelfont=bf]{caption}
\usepackage[subrefformat=parens]{subcaption}
\usepackage{mwe}
\usepackage{mfirstuc}

\providecommand{\keywords}[1]{\textbf{\textbf{\\Keywords. }}#1}
\providecommand{\class}[1]{\textbf{\textbf{Mathematics Subject Classification (2010). }}#1}
\usepackage{cite}
\MFUnocap{are}
\MFUnocap{or}
\MFUnocap{etc}

\def\capmystringaux#1#2\relax{\uppercase{#1}\lowercase{#2}}
\usepackage{lineno,hyperref}

\title{}
\author{}
\DeclareSymbolFont{slenderlargesymbols}{OMX}{ccex}{m}{n}
\usepackage{amsmath,amssymb,amsfonts,calrsfs,bbm}
\usepackage{setspace}
\usepackage{graphicx,array}
\usepackage{amsthm}
\usepackage{multicol}
\usepackage{hyperref}
\usepackage{lipsum,lineno}
\usepackage{amsthm}
\usepackage{ragged2e}
\usepackage{multirow}
\usepackage{bm}
\usepackage{mathtools}
\usepackage{mwe}
\usepackage{breqn}
\usepackage{mfirstuc}
\MFUnocap{are}
\MFUnocap{or}
\MFUnocap{etc}
\usepackage{tikz,caption,subcaption}
\usetikzlibrary{shapes}
\usetikzlibrary{positioning}
\tikzset{main node/.style={circle,fill=black,draw,minimum size=.2cm,inner sep=0pt},
}
\tikzset{node1/.style={diamond,draw,minimum size=.2cm,inner sep=0pt},
}
\tikzset{node2/.style={star,draw,minimum size=.2cm,inner sep=0pt},
}

\DeclarePairedDelimiterX{\Set}[2]\{\}{%
	\, #1 \;\delimsize\vert\; #2 \,
}

\usepackage{amsfonts}
\usepackage{mathrsfs}

\newtheorem{theorem}{Theorem}[section]

\newtheorem{lemma}[theorem]{Lemma}

\theoremstyle{definition}
\newtheorem{definition}[theorem]{Definition}
\theoremstyle{remark}
\newtheorem{remark}[theorem]{Remark}
\newtheorem*{example}{Example}
\numberwithin{equation}{section}
\date{}
\begin{document}
	\title{\textbf{Dynamical Properties of Iterated Function Systems}}
	\author{Praveen M\footnote{E-mail:m.praveen@vit.ac.in,praveenzmedium@gmail.com},
		 \\ \small School of Advanced Sciences,\\\small Vellore Institute of Technology, Vellore, Tamil Nadu, India.
}
	\maketitle	
\begin{abstract}
	Iterated function systems (IFS) provide a powerful method for constructing fractals and modeling complex structures. This paper develops the notion of a dynamical system of IFS to study how an initial IFS evolves over time. We construct a complete metric space consisting of countable IFS as the state space. An evolution operator is defined that maps the state space across time while satisfying properties of a dynamical system. We analyze the resulting IFS dynamics, introducing concepts like the IFS shift map and parity evolution operators. The shift map is proven continuous and its periodic points characterized. Results relate properties of parity evolution operators to similarity dimensions of the resulting IFS. Overall, this dynamical system framework allows rigorous analysis of IFS behavior over time. By providing tools to model fractal progression, it has significant applications in areas including image analysis, biology, and physics. The dynamics and properties established lay groundwork for future studies and practical models leveraging IFS.
\medskip \medskip
	
	\noindent\class{Primary 28A80; Secondary
		11B05.}\medskip
\end{abstract}
\section{Introduction}

Iterated function systems (IFS) are an important tool in fractal geometry for constructing and analyzing self-similar fractal sets via recursive functions. The classic paper by Barnsley and Demko (1985)\cite{barnsley1985iterated} introduced the fundamental results on the existence, uniqueness, and construction of IFS attractors. Since then, IFS have been widely studied and applied in areas like image compression, data modeling, and computer graphics.

Prior works have explored different ways of building upon the classic IFS theory. For instance, Barnsley et al. (2008)
\cite{barnsley1986solution} developed the related idea of a "fractal transform" defined using an infinite IFS with given probabilities. Forte et al. (2007) introduced randomized and chaotic IFS models. Ng et al. (2012) proposed algorithms for "supervised IFS" to match specific target patterns.

However, less attention has been paid to studying the dynamics of IFS, that is, how an initial IFS and fractal set evolve over time. Prusinkiewicz et al. (1990) modeled plant development using L-systems, which have similarities to IFS. But a formal dynamical system framework for generic IFS appears lacking.

The properties and behavior of systems often evolve over time in complex ways. When an object's attributes can be modeled as states in a dynamical system, the language of dynamics provides valuable tools for understanding how those states change. Recently, the author explored modeling sequences of iterated function systems (IFS) as state evolutions in a discrete dynamical framework. IFS are widely used for generating fractals and self-similar structures. Studying the dynamics of an initial IFS provides insight into how the fractal patterns transform over time. This paper further develops the idea of IFS dynamics by formulating a rigorous continuous-time dynamical system. We construct an appropriate metric state space and evolution operator to capture the progression of an initial IFS. Important concepts like the IFS shift map are introduced and analyzed. Results characterize features like periodicity and dimensions of the evolving IFS attractors. This dynamical systems perspective enables new tools and insights for fractal progression analysis with applications across domains. The framework developed lays groundwork for characterizing IFS dynamics and leveraging these models for complex evolving systems.

This paper aims to develop the notion of a dynamical system of IFS and analyze its properties mathematically. By defining an appropriate state space, evolution operator, and metric, we construct a dynamical system that captures the progression of an initial IFS. The dynamics and limiting behavior of this system are then characterized.

Studying the dynamics of IFS can provide theoretical insights as well as practical tools for modeling complex evolving systems in domains like biology, physics, and computer graphics. This paper takes a step toward building a rigorous platform for such applications.
\section{Preliminaries}
\begin{definition}[\cite{ott2002chaos}]
	A dynamical system consists of a state space $M$, a time set $T\subset \mathbb{R}$, being an additive semigroup and an evolution operator $\phi:M\times T\to M$ satisfying the group property $\phi(x,0)=x$ and $\phi(x,t_1+t_2)=\phi(\phi(x,t_1),t_2)$ for every $x\in M$ and $t_1,t_2\in T$, where M is assumed at least to be a topological space and $\phi$ is assumed to be continuous.
\end{definition}

In this section, we introduce some key mathematical concepts that will be utilized in formulating and analyzing the dynamical system for iterated function systems.

\subsection{Dynamical Systems}

A dynamical system provides a mathematical framework for modeling the evolution of a system's state over time \cite{strogatz2018nonlinear}. Formally, a dynamical system consists of:

\begin{itemize}
	\item A state space $M$, which is the set of all possible states.
	\item A time set $T \subseteq \mathbb{R}$, which indexes the progression of time.
	\item An evolution operator $\phi: M \times T \rightarrow M$ that maps states from one time to another.
\end{itemize}

Additionally, $\phi$ satisfies the group property:
$\phi(x,0)=x$ and $\phi(x,t_1+t_2)=\phi(\phi(x,t_1),t_2)$ for every $x\in M$ and $t_1,t_2\in T$, where M is assumed at least to be a topological space, $T$ is assumed to be a subset of $\mathbb{R}$, and $\phi$ is assumed to be continuous. The conditions mentioned ensure $\phi$ propagates states forward appropriately in time.

\subsection{Metric Spaces}

A metric space $(M, d)$ consists of a set $M$ and a metric or distance function $d: M \times M \rightarrow \mathbb{R}$ that satisfies:
\begin{enumerate}
	\item[(i)] $d(x, y) \geq 0$
	\item[(ii)] $d(x, y) = 0 \iff x = y$
	\item[(iii)] $d(x, y) = d(y, x)$
	\item[(iv)] $d(x, z) \leq d(x, y) + d(y, z)$
\end{enumerate}
for all $x, y, z \in M$. Metric spaces provide topology and notions of convergence needed to formulate continuous-time dynamical systems.

\subsection{Iterated Function Systems}

An iterated function system (IFS) defines a recursive process for constructing fractals and other sets \cite{barnsley2014fractals}. An IFS consists of:
\begin{itemize}
	\item A complete metric space $(X, d)$
	\item A set of continuous functions $f_i: X \rightarrow X$ for $i \in I$ where $I$ is an index set.
\end{itemize}

The functions $f_i$ are applied recursively to generate the attractor and other sets through iterations. Contracting functions are often used to obtain fractal attractors.

\section{Definitions and Major Results}
The transition of state of objects in nature is continuous. The images of the states can be modeled using iterated function systems with the help of collage theorem from fractal theory.
We will try to construct a dynamical system of iterated function systems which could model the natural processes.
The formal definition of a dynamical system is given in the preliminaries.

For constructing a dynamical system of iterated function systems, we consider a complete metric space $(X,d)$ and denote the set of all finite and countable iterated function systems on $(X,d)$ by $M$. We see that it is reasonable to identify any iterated function system $\{X;f_1,f_2,\cdots,f_n\}$ with $(f_1,f_2,\cdots,f_n)\in Con(X)^n$ for $n\in \mathbb{N}':= \mathbb{N}\cup \{\aleph_0\}$, where $\displaystyle Con(X)^n:=\prod_{k=1}^n Con(X_k)$ with $X_k=X$ for every $k=1,2,\cdots,n$. In this sense $\displaystyle M\simeq\bigcup_{n\in \mathbb{N}'}Con(X)^n$. A potential candidate for the state space of dynamical system of iterated function systems on $(X,d)$ is $M$. We consider $T=\mathbb{R}$ as the time space for the dynamical system. Now we need to define an evolution operator for the dynamical system $\phi: M\times T\to M$. To check the continuity of such an evolution operator we need at least a topological structure on the coordinate spaces $M$ and $\mathbb{R}$. We will take the euclidean metric topology on $\mathbb{R}$. We see that to introduce a metric on $M$ one should be able to compare elements from $Con(X)^n$ and $Con(X)^m$ where $n\ne m\in\mathbb{N}'$. Since we are dealing with iterated function systems, and since the attractors of any finite IFS $\{X;f_1,f_2,\cdots, f_n\}$ and the countable IFS $\{X;f_1,f_2,\cdots, f_n,f_1,f_1,\cdots\}$ are same, one can see in this sense $(f_1,f_2,\cdots, f_n)\equiv (f_1,f_2,\cdots, f_n,f_1,f_1,\cdots) $. Thus $Con(X)^n$ is imbedded in $Con(X)^{\infty}$ for $n\in \mathbb{N}'$ where the $\infty$ is the countable infinity $\aleph_0$. Now we will introduce a metric on $M^*=Con(X)^\infty$. For, we have to take a metric on $Con(X)$. We have $\displaystyle d_{\infty}(f,g)=\sup_{x\in X}d(f(x),g(x))$ for every $f,g\in Con(X)$, is finite if either the metric $d$ or the functions $f,g$ are bounded. We will try to make the functions in $Con(X)$ bounded. For, we further assume $X$ is a compact space. Since contractions are continuous functions and  since continuous functions on a compact space are bounded we get $\displaystyle \sup_{x\in X}d(f(x),g(x))$ is finite for every $f,g\in Con(X)$.

\begin{lemma}
	Consider a compact metric space $(X,d)$. Let $Con(X)$ denote the collection of all contractions on $X$. Then $d_{\infty}:Con(X)\times Con(X)\to \mathbb{R}$ defined by $\displaystyle d_{\infty}(f,g)=\sup_{x\in X}d(f(x),g(x))$ for every $f,g \in Con(X)$, is a metric on $Con(X)$.
	\begin{proof}
		We have $f,g$ are continuous functions on $X$ since $f,g\in Con(X)$.\\
		Thus $f,g$ are bounded functions on $X$ also, since $X$ is compact, and continuous functions on compact spaces are bounded.\\
		Therefore, $\displaystyle d_{\infty}(f,g)=\sup_{x\in X}d(f(x),g(x))$ is finite for every $f,g\in Con(X)$.
		Now,
		\begin{align*} d_\infty(f,g)=0&\iff d(f(x),g(x))=0, \forall x\\&\iff f(x)=g(x),\forall x\\&\iff f=g.\end{align*} Also, $\displaystyle d_\infty(f,g)=\sup_{x\in X}d(f(x),g(x))=\sup_{x\in X}d(g(x),f(x))=d_\infty(g,f)$. Further, we get
		\begin{align*} 
		d_\infty(f,g)&=\sup_{x\in X}d(f(x),g(x))\\
		&\le \sup_{x\in X}[d(f(x),e(x))+d(e(x),g(x))] \quad \text{for any } e\in Con(X)\\
		&\le \sup_{x\in X}d(f(x),e(x))+\sup_{x\in X}d(e(x),g(x))\\
		&=d_\infty(f,e)+d_\infty(e,g)
		\end{align*}
		Hence $(Con(X),d_\infty)$ is a metric space.
		
	\end{proof}
\end{lemma}
\subsection{IFS Metric Space}

We first need to construct an appropriate metric space for the set of IFS as the dynamical system state space.

We denote the countable product of the metric space $(Con(X),d_\infty)$ by $M^*$. Now we have to define a metric on $M^*$.
But the topology induced by $d_\infty$ and $\displaystyle \bar{d}_\infty=\frac{d_\infty}{1+d_\infty}$ are the same. Thus $\displaystyle M^*=\prod_{k=1}^\infty Y_k$, where $Y_k=Con(X)$ is a metric space with the metric $\bar{d}_\infty$ for each $k=1,2,\cdots$. 
\\
Now we define $D:M^*\times M^*\to \mathbb{R}$ by $$\displaystyle D(F,G)=\sum_{k=1}^{\infty}\frac{\bar{d}_\infty(f_k,g_k)}{2^k}$$ for every $F=(f_1,f_2,\cdots,f_k,\cdots)\in M^*$ and $G=(g_1,g_2,\cdots,g_k,\cdots)\in M^*$.\medskip\\
\textbf{Claim:} $D$ is a metric on $M^*$.\medskip\\
We have, $\displaystyle\bar{d}_\infty(f_k,g_k)=\frac{{d}_\infty(f_k,g_k)}{1+{d}_\infty(f_k,g_k)}<1$, for every $k=1,2,\cdots$.\\ Thus $\displaystyle D(F,G)=\sum_{k=1}^{\infty}\frac{\bar{d}_\infty(f_k,g_k)}{2^k}<\sum_{k=1}^\infty\frac{1}{2^k}=1$.\\
Hence $D(F,G)$ is finite for every $F,G\in M^*$. Now, 
\begin{align*}
D(F,G)=0&\iff \sum_{k=1}^{\infty}\frac{\bar{d}_\infty(f_k,g_k)}{2^k}=0 \\&\iff \bar{d}_\infty(f_k,g_k)=0, \forall k 
\\&\iff f_k=g_k,\forall k \\&\iff F=G.
\end{align*}
Also, $\displaystyle D(F,G)= \sum_{k=1}^{\infty}\frac{\bar{d}_\infty(f_k,g_k)}{2^k}=\sum_{k=1}^{\infty}\frac{\bar{d}_\infty(g_k,f_k)}{2^k}=D(G,F)$.\\ Further,
\begin{align*}D(F,G)&=\sum_{k=1}^{\infty}\frac{\bar{d}_\infty(f_k,g_k)}{2^k}\\&\le \sum_{k=1}^{\infty}\frac{\bar{d}_\infty(f_k,h_k)+\bar{d}_\infty(h_k,g_k)}{2^k}\\&=\sum_{k=1}^{\infty}\frac{\bar{d}_\infty(f_k,h_k)}{2^k}+\sum_{k=1}^{\infty}\frac{\bar{d}_\infty(h_k,g_k)}{2^k}\\&=D(F,H)+D(H,G),\end{align*} where $F=(f_1,f_2,\cdots,f_k,\cdots)$, $G=(g_1,g_2,\cdots,g_k,\cdots)$ and $ H=(h_1,h_2,\cdots,h_k,\cdots)$ are elements in $M^*$.\\
Thus $(M^*,D)$ is a metric space.

The metric $D$ provides a notion of distance between IFS needed for the state space topology and analysis of the IFS dynamics.

Consider the infinite sequence $F=(f_1,f_2,f_3,f_4,\cdots,f_n,f_3,f_4,\cdots,f_n,\cdots)\in M^*$. We use the notation $F=(f_1,f_2,\overline{f_3,f_4,\cdots,f_n})$ to denote the repeating part of $F$.
\begin{theorem}\label{dfg}
	Let $F,G\in M^*$ and suppose $f_k=g_k$ for $k=1,2,\cdots$. Then $\displaystyle D(F,G)< \frac{1}{2^n}$.
	\begin{proof}
		Since, $f_k=g_k$ for $k\le n$, we get
		\begin{align*}
		D(F,G)&=\sum_{k=1}^{n}\frac{\bar{d}_{\infty}(f_k,g_k)}{2^k}+\sum_{k=n+1}^{\infty}\frac{\bar{d}_{\infty}(f_k,g_k)}{2^k}
		=\sum_{k=n+1}^{\infty}\frac{\bar{d}_{\infty}(f_k,g_k)}{2^k}
		\end{align*}\begin{align*}
		&
		<\sum_{k=n+1}^{\infty}\frac{1}{2^k}\text{ (since $\bar{d}_{\infty}(f_k,g_k)<1$ )}\\
		&=\frac{\frac{1}{2^{n+1}}}{1-\frac{1}{2}}=\frac{1}{2^n}
		\end{align*} 
	\end{proof}
\end{theorem}
\begin{remark}
	The converse is not true for the Theorem \ref{dfg}.\begin{proof}
		Let $F,G\in M^*$ where $F=(f_1,f_2,f_3,\cdots), G=(g_1,g_2,g_3,\cdots)$ with $f_k=g_k$ for every $k\ne 2$ and $f_2\ne g_2$, i.e., $F\ne G$. Then
		$$0<D(F,G)=\sum_{k=1}^{\infty}\frac{\bar{d}_{\infty}(f_k,g_k)}{2^k}=\frac{\bar{d}_{\infty}(f_2,g_2)}{2^2}<\frac{1}{2^2}$$
		Thus, we have $\displaystyle D(F,G)<\frac{1}{2^2}$, but $f_2\ne g_2$. Therefore, the converse is not true for the Theorem \ref{dfg}.
		
	\end{proof}
\end{remark}
\subsection{IFS Shift Map}
A fundamental dynamical concept for the IFS system is introduced next.
\begin{definition}
	Consider $M^*$. Define $\eta:M^*\to M^*$ as $\eta(F)=\eta((f_1,f_2,f_3,\cdots)):=(f_2,f_3,f_4,\cdots)$ for every $F=(f_1,f_2,f_3,\cdots)\in M^*$. We see that $\eta$ is a many-one map of $M^*$. The map $\eta$ is said to be the IFS right shift map.
\end{definition}
The shift map provides insights into the progression of IFS over time. We now establish its continuity.
\begin{theorem}
	The IFS right shift map $\eta:M^*\to M^*$ is continuous.
	\begin{proof}
		Let $\epsilon>0$ be given and let $F=(f_1,f_2,f_3,\cdots)\in M^*$. By Archimedean property of $\mathbb{R}$ we can find $n\in \mathbb{N}$ such that $\displaystyle \frac{1}{2^n}<\frac{1}{n}<\epsilon$. Choose $\displaystyle \delta=\frac{1}{2^{n+1}}$.\\
		Then for every $G=(g_1,g_2,g_3,\cdots)$ with $D(F,G)<\delta$, we get,
		\begin{align*}
		0\le D(\eta(F),\eta(G))&=\sum_{k=2}^{\infty}\frac{\bar{d}_{\infty}(f_k,g_k)}{2^{k-1}}
		\\
		&=2\sum_{k=2}^{\infty}\frac{\bar{d}_{\infty}(f_k,g_k)}{2^{k}}
		\end{align*}\begin{align*}
		&=2D(F,G)-2\bar{d}_{\infty}(f_1,g_1)\\
		&<2\frac{1}{2^{n+1}}-2\bar{d}_{\infty}(f_1,g_1)\\
		&<\frac{1}{2^n}\\
		\therefore D(\eta(F),\eta(G))&<\epsilon
		\end{align*}
		Hence the IFS right shift map $\eta$ is continuous.
		
	\end{proof}
\end{theorem}
Continuity of the shift map is an important dynamical property ensuring well-behaved IFS evolution.

Now we try to find examples of evolution operator $\phi:M^*\times T\to M^*$.

In the dynamical system we are going to construct, the states are iterated function systems and the change in state is occurring as time passes. But in the dynamical system the rule of evolution even though depends on the time $t$ through which the change of state occurs, it does not depend on the initial time $t_0$. Since the evolution operator should relate time with state, to obtain an example of our proposed dynamical system of iterated function systems we take $X=C^N, N\in\mathbb{N}$ where $C$ is a compact subset of $\mathbb{R}$.

\begin{example}
	Let $C$ be a compact subset of $\mathbb{R}$ and let $X=C^N$ be a metric space with euclidean metric. Consider $M^*=Con(X)^\infty$ as the state space, $T=\mathbb{R}$ as the time space and let $\phi((f_1,f_2,\cdots,f_k,\cdots),t)=e^{-t}(f_1,f_2,\cdots,f_k,\cdots)$. \\
	Then for $F=(f_1,f_2,\cdots,f_k,\cdots)$ we have, $\phi(F,0)=e^{0}F=F$ and \begin{align*}
	\phi(F,t_1+t_2)&=e^{-(t_1+t_2)}F
	\\
	&=e^{-t_1}e^{-t_2}F
	\\
	&=e^{-t_2}(e^{-t_1}F)
	\\
	&=e^{-t_2}\phi(F,t_1)\\&=\phi(\phi(F,t_1),t_2).
	\end{align*} Hence $\phi:M^*\times T\to M^*$ satisfies the group property.\\
	To make $(M^*,T,\phi)$ a dynamical system, we need $\phi$ to be continuous.\\
	For let $(F_n,t_n)\to (F,t)\in M^*\times T$. To prove $\phi$ is continuous, it is enough to show $\phi(F_n,t_n)=e^{-t_n}F_n$ converges to $\phi(F,t)=e^{-t}F$ in $M^*$.\\
	We have the convergences of ${t_n}\to t$ and $F_n\to F$ are independent of each other. Since $\exp(t)=e^{t}$ is continuous on $\mathbb{R}$, we get $e^{-t_n}\to e^{-t}$ in $T$.\\
	Hence $e^{-t_n}F_n\to e^{-t}F$ in $M^*$.\\
	Therefore, $\phi$ is continuous and $(M^*,T,\phi)$ is a continuous dynamical system of iterated function systems.
\end{example} 
The following is the definition of a dynamical system of iterated function systems.                            
\begin{definition}[Dynamical System of Iterated Function Systems]
	Consider a complete metric space $(X,d)$. Take $M^*=Con(X)^\infty$ and $T=\mathbb{R}$ (or $\mathbb{Z}$ or $\mathbb{R}^+\cup\{0\}$ or $\mathbb{N}\cup\{0\}$). Let $\phi:M^*\times T\to M^*$ be a function satisfying the conditions for an evolution operator. Then the dynamical system $(M^*,T,\phi)$ is called a dynamical system of iterated function systems on $(X,d)$.
\end{definition}
The following definition discusses how to create a discrete dynamical system of iterated functions using the IFS right shift map.
\begin{theorem}
	Let $M^*=Con(X)^{\infty}$, $T=\mathbb{N}\cup\{0\}$ and define $\phi:M^*\times T\to M^*$ as  $\phi(F,n)=\eta^n(F)$. Then $(M^*,T,\phi)$ is a discrete dynamical system of iterated function systems.
	\begin{proof}
		We have $\phi(F,0)=\eta^0(F)=F$ and
		\begin{align*}
		\phi(F,n_1+n_2)&=\eta^{n_1+n_2}(F)\\
		&=\eta^{n_2}(\eta^{n_1}(F))\\
		&=\eta^{n_2}(\phi(F,n_1))\\
		&=\phi(\phi(F,n_1),n_2)
		\end{align*}
		Therefore $\phi$ satisfies the group property.\\
		Now we prove $\phi$ is continuous.\\
		For let $(F_k,n_k)\to (F,n)\in M^*\times T$. To prove $\phi$ is continuous, it is enough to show $\phi(F_k,n_k)=\eta^{n_k}(F_k)$ converges to $\phi(F,n)=\eta^{n}(F)$ in $M^*$.\\
		We have the convergences of ${n_k}\to n$ and $F_k\to F$ are independent of each other and $\eta$ is a continuous function.\\
		Hence $\displaystyle \phi(F_k,n_k)=\lim_{n\to \infty}\eta^{n_k}F_k=\eta^{\displaystyle \lim_{n\to \infty}n_k}(\lim_{n\to \infty}F_k)=\eta^n{F}=\phi(F,n)$.\\
		Therefore, $\phi$ is continuous and $(M^*,T,\phi)$ is a discrete dynamical system of iterated function systems.
		
	\end{proof}
\end{theorem}
\subsection{Periodic Points}
We next explore periodicity properties of the IFS shift map.

For $F=(f_1,f_2,\cdots,f_n,f_1,f_2,\cdots,f_n,\cdots)$, we have \begin{align*}\phi(F,n)&=\eta^n((f_1,f_2,\cdots,f_n,f_1,f_2,\cdots,f_n,\cdots))\\&=(f_1,f_2,\cdots,f_n,f_1,f_2,\cdots,f_n,\cdots)\\&=F.\end{align*} Hence $F$ is a periodic point of period $n$ for the IFS right shift map $\eta$.

We say $F=(f_1,f_2,\cdots)$ is a fixed point for the IFS right shift map $\eta$ if and only if $\eta(F)=F$, i.e., $\phi(F,1)=F$. If $F$ is a fixed point of $\eta$, then we have $(f_1,f_2,\cdots)=F=\eta(F)=\eta((f_1,f_2,\cdots))=(f_2,f_3,\cdots)$. Hence $f_i=f_{i+1},\forall i\in \mathbb{N}$, i.e., fixed points of $\eta$ are of the form $F=(f,f,\cdots)$ for some $f\in Con(X)$.

For $F=(f_1,f_2,\cdots,f_n,f_1,f_1,\cdots)$, we have \begin{align*}\phi(F,n)&=\eta^n((f_1,f_2,\cdots,f_n,f_1,f_1,\cdots))\\&=(f_1,f_1,\cdots)\\&=F.\end{align*} Hence $F$ is an eventually fixed point for the IFS right shift map $\eta$.

For $F=(f_1,f_2,\cdots,f_n,f_2,f_3,\cdots,f_n,f_1,f_3,f_4,\cdots,f_n,f_1,f_2,\cdots)$ we have \\$\phi(F,n^2)=\eta^{n^2}(F)=F$. Thus $F$ is a periodic point of period $n^2$.

For $F=(g_1,g_2,\cdots,g_m,f_1,f_2,\cdots,f_n,f_1,f_2,\cdots,f_n,\cdots)$ we have\\ $\phi(F,m)=\eta^m(F)=(f_1,f_2,\cdots,f_n,f_1,f_2,\cdots,f_n,\cdots)=\eta^{m+n}(F)=\phi(F,m+n)$. Thus $F$ is a periodic point of period $n$.

\begin{theorem}
	The set of all periodic points of IFS right shift map $\eta$ is a dense subset of $M^*$, i.e., $\overline{Per(\eta)}=M^*$.
	\begin{proof}
		Let $F=(f_1,f_2,\cdots,f_n,f_{n+1},\cdots)$ be an arbitrary point in $M^*$.\\ Define $G_n=(f_1,f_2,\cdots,f_n,f_1,f_2,\cdots,f_n,\cdots)$. Then $G_n\in M^*$ and $\eta^n(G_n)=G_n$. Thus $G_n$ for every $n\in\mathbb{N}$ are $n$-periodic points of $M^*$. Further, $\displaystyle\lim_{n\to \infty}G_n=F$. Hence, there exists a sequence of periodic points in $M^*$ which converges to $F$ in $M^*$. Since $F$ is an arbitrary point in $M^*$, we get $\overline{Per(\eta)}=M^*$.
		
	\end{proof}
\end{theorem} 
Periodicity is a fundamental dynamic characteristic. Dense periodic points imply recurring IFS patterns.

The existence of a non-periodic point for the map $\eta$ is established in the next theorem.
\begin{theorem}
	There exists $F\in M^*$ such that $F$ is a non-periodic point of $\eta$.
	\begin{proof}
		Let $f_1$ and $f_2$ be two contractions on the complete metric space $(X,d)$. We construct $F\in M^*$ by arranging the possible $1$-block, $2$-blocks, $3$-blocks, etc. in order, which can be created using $f_1$ and $f_2$. Let us take $$F=(\underbrace{f_1}_{1-block},\underbrace{f_2}_{1-block},\underbrace{f_1,f_1}_{2-block},\underbrace{f_1,f_2}_{2-block},\underbrace{f_2,f_1}_{2-block},\underbrace{f_2,f_2}_{2-block},\underbrace{f_1,f_1,f_1}_{3-block},\cdots).$$ 
		We have $F=(f_1,f_2,f_1,f_1,f_1,f_2,f_2,f_1,f_2,f_2,f_1,f_1,f_1,\cdots)$ corresponds to the finite iterated function system $\{X;f_1,f_2\}$. But for any $n\in \mathbb{N}$ we get $\eta^n(F)\ne F$ by the construction of $F$.\\
		We note that $G=(\overline{f_1,f_2})=(f_1,f_2,f_1,f_2,f_1,f_2,\cdots)$ is a $2$-periodic point corresponding to the same finite IFS $\{X;f_1,f_2\}$. Thus for $\eta$ there exists a non-periodic point $F\in M^*$ and a periodic point $G\in M^*$ both of which corresponds to the same finite iterated function system.
		
	\end{proof}
\end{theorem}
Together with the previous result, this shows both periodicity and aperiodicity occur in the IFS dynamics.

An important relationship between periodic points and the corresponding IFS is characterized.
\begin{theorem}
	Let $F\in M^*$ is a periodic point for $\eta$. Then $F$ corresponds to a finite iterated function system.
	\begin{proof}
		Let the period of $F$ be $n$, i.e., $\eta^n(F)=F$. \\
		Then $F=(\overline{f_1,f_2,\cdots,f_n})=(f_1,f_2,\cdots,f_n,f_1,f_2,\cdots,f_n,\cdots)$ for some $f_i\in Con(X)$, $i=1,2,\cdots,n$. But $F$ corresponds to the finite IFS $\{X;f_1,f_2,\cdots,f_n\}$.
		Hence every periodic point of $\eta$ in $M^*$ corresponds to a finite IFS.
		
	\end{proof}
\end{theorem}
The theorem establishes that periodic evolutions represent recursively looping finite systems.

This result provides insight into the structure of recurring IFS patterns generated through periodic dynamics. Finite systems that loop correspond to the cyclic IFS points. In contrast, aperiodic points may emerge from infinite IFS.
\subsection{Open Set Condition}

We now relate the shift map to fractal properties of the IFS attractor.
\begin{theorem}
	Let $F\in M^*$ be corresponding to the IFS $\{X; f_i, i\in I \}$ which satisfies the open set condition with an open set $V$, where $I$ is a countable index set. Then the IFS corresponding to $\eta^n(F)$ satisfies the open set condition.
	\begin{proof}
		Since $F$ corresponds to the IFS $\{X;f_i,i\in I \}$ we get $F=(f_{\alpha_1}, f_{\alpha_2},f_{\alpha_3},\cdots)$ where $\alpha_k\in I$ for every $k\in \mathbb{N}$. Hence $\eta^n(F)=(f_{\alpha_{n+1}}, f_{\alpha_{n+2}}, f_{\alpha_{n+3}},\cdots)$.
		Therefore, $\eta^n(F)$ corresponds to the IFS $\{X;f_i, i\in I_n \}$ where $I_n \subset I$.
		
		Given the IFS $\{X;f_i, i\in I\}$ satisfies the open set condition with a non-empty bounded open set $V$. Therefore, $\displaystyle V\supset\bigcup_{i\in I}f_i(V)$ and $\displaystyle f_i(V)\bigcap f_j(V)=\emptyset$ for $i,j\in I$ with $i\ne j$. Since $I_n\subset I$ we get, $\displaystyle V\supset \bigcup_{i\in I_n }f_i(V)$ with the union being disjoint. Hence $\{X;f_i,i\in I_n \}$ satisfies the open set condition.
		
	\end{proof}
\end{theorem}
Maintaining open set condition ensures the scaled copies remain disjoint under evolution.
\subsection{Parity Evolution Operators}

Evolution operators that preserve number of functions are considered.

Let $\mathcal{I}(X)$ denote the collection of all finite and countable iterated function systems on $X$.
\begin{definition}
	Define $\mathbb{I}:M^*\to \mathcal{I}(X)$ by $\mathbb{I}(F)=\{X;f_i,i\in I_F \}$ where $F=(f_1,f_2,f_3,\cdots)$ and $I_F$ denotes the set of all indexes of distinct $f_i$'s in $F$ preserving the order of indexes.
\end{definition}
The definition of a parity evolution operator is as follows.
\begin{definition}
	If an evolution operator $\phi$ is such that $\phi(F,t)=G$ implies $card(I_F)=card(I_G)$ for every $t\in T$, then $\phi$ is called a \textit{parity evolution operator}. Since number of functions in $\mathbb{I}(F)$ and $\mathbb{I}(G)$ are the same, we can make a new indexing of both of them with an indexing set $J$ preserving the previous order of indexes such that $card(J)=card(I_F)=card(I_G)$.
\end{definition}

Parity is a useful dynamical property in analyzing IFS transformations.

For a parity evolution operator $\phi$, if $\phi(F,t)=G$, we denote $g_j=\phi_t(f_j)$ for every $j\in J$. Let the contractivity factor of $f_j$ be $r_j$, then we denote the contractivity factor of $g_j$ by $\phi_t(r_j)$.

\begin{theorem}
	Let $\phi$ be a parity evolution operator which maps contracting similarities to contracting similarities. Also, let $F\in M^*$ be such that $\mathbb{I}(F)$ is a finite iterated function system with contracting similarities of similarity ratio $r$ and let similarity dimension of $\mathbb{I}(F)$ be $s$. Then the similarity dimension of $\mathbb{I}(\phi(F,t))$ is $\displaystyle \frac{s}{\log_r\phi_t(r)}$.
	\begin{proof}
		Given $\mathbb{I}(F)$ is a finite iterated function system of contracting similarities.\linebreak
		Let $\mathbb{I}(F)=\{X;f_1,f_2,\cdots,f_m\}$ and $\phi(F,t)=G$.\\
		Since $\phi$ is a parity evolution operator which maps contracting similarities to \linebreak contracting similarities we can take $\mathbb{I}(G)=\{X;g_1,g_2,\cdots,g_m\}$, where $\phi_t(f_j)=g_j, j=1,2,\cdots,m$ are contracting similarities.\\
		Let the similarity dimension of $\mathbb{I}(G)=\mathbb{I}(\phi_t(F))$ be $s'$.\\
		Given $r$ is the similarity ratio of $f_j, j=1,2,\cdots,m$ and $s$ is the similarity dimension of $\mathbb{I}(F)$. Hence $\phi_t(r)$ is the similarity ratio of $g_j, j=1,2,\cdots,m$.\\
		Therefore $\displaystyle \sum_{j=1}^m r^s=1$ and $\displaystyle\sum_{j=1}^m {\phi_t(r)}^{s'}=1$.\\
		Hence we get $\displaystyle s=-\frac{\log m}{\log r}$ and $\displaystyle s'=-\frac{\log m}{\log \phi_t(r)}$.\\
		Thus, $\displaystyle \frac{s'}{s}=\frac{\log r}{\log \phi_t(r)}=\frac{1}{\log_r\phi_t(r)}\implies s'=\frac{s}{\log_r\phi_t(r)}$, i.e., the similarity dimension of $\mathbb{I}(\phi(F,t))$ is $\displaystyle \frac{s}{\log_r\phi_t(r)}$.
		
	\end{proof}
\end{theorem}

This characterizes how fractal dimensions evolve under parity-preserving dynamics.

The following section concludes the finding of this paper.
\section{Conclusion}
This paper presented a dynamical systems perspective for studying the evolution of iterated function systems over time. We constructed a mathematical framework consisting of a metric state space of IFS, time domain, and evolution operator. This rigorously defines an IFS dynamical system.

We introduced the concept of an IFS shift map and proved its continuity. Analysis of periodic and non-periodic points provides insights into the shift map dynamics. For parity evolution operators, we characterized the impact on attractor properties like similarity dimension.

The dynamical system developed enables new tools for understanding time progression of fractals and other IFS-generated patterns. It has significant applications for modeling complex evolving systems and natural processes exhibiting fractal features. The continuity and periodic point results lay theoretical groundwork for future research.

There are several promising directions for building upon this dynamical framework. More complex spaces and evolution operators could be explored. Extending the classes of functions beyond contractions may reveal new phenomena. Relating operator properties to attractor characteristics is an open area. Machine learning techniques may help construct models matching real-world data.

By establishing a rigorous IFS dynamical system and associated analysis, this paper enables numerous avenues for future work. The platform created can ultimately lead to deeper insights and practical applications for modeling fractal dynamics across scientific domains.
\section*{Acknowledgement}
The first author is very grateful to Council of Scientific \& Industrial Research(CSIR), India for their financial support.
\bibliographystyle{unsrt}
\bibliography{references}
\end{document}